\documentclass{amsart}

\usepackage{amsthm}
\usepackage{amssymb}
\usepackage{amsmath}
\usepackage{amsfonts}
\usepackage{amscd}
\usepackage{enumerate}

\title[Separating twists and the Magnus representation]
{Separating twists and the Magnus representation\\ of the Torelli group} 

\author{Thomas Church}
\address{Department of Mathematics\\
5734 S. University Ave.\\
Chicago, IL 60637}
\email{tchurch@math.uchicago.edu\newline\indent http://math.uchicago.edu/$\sim$tchurch/}

\author{Aaron Pixton}
\address{Mathematics Department, Princeton University\\
Princeton, NJ 08544}
\email{apixton@math.princeton.edu}

\theoremstyle{plain}
\newtheorem{conjecture}{Conjecture}
\newtheorem{question}[conjecture]{Question}
\newtheorem{theorem}{Theorem}[section]
\newtheorem{proposition}[theorem]{Proposition}
\newtheorem{lemma}[theorem]{Lemma}
\newtheorem{corollary}[theorem]{Corollary}

\newtheorem*{freelem}{Lemma~5.1}
\newtheorem*{mainthm}{Theorem~6.1}

\theoremstyle{definition}
\newtheorem{definition}[theorem]{Definition}
\newtheorem*{rmk}{Remark}

\newcommand{\Z}{\mathbb{Z}}

\newcommand{\Q}{\mathbb{Q}}

\newcommand{\freetext}{Let $V$ be a finite-dimensional vector space
  over a field $K$ of characteristic zero.  Let $A, B\in\End(V)$
  satisfy $A^2 = B^2 = 0$, and suppose that $\tr(AB)\in K$ is
  transcendental over $\Q$. Then the endomorphisms $1+A$ and $1+B$ generate a free
  subgroup of rank $2$ in $\Aut(V)$.}
\newcommand{\tfaetext}{Let $T_C
  = T_{\gamma_1}^{n_1}\cdots T_{\gamma_k}^{n_k}$ and $T_D =
  T_{\delta_1}^{m_1}\cdots T_{\delta_l}^{m_l}$ be positive separating
  multitwists. Let $c_1,\ldots,c_k$ and $d_1,\ldots,d_l\in\Hcover$ be
  lifts of $\gamma_1,\ldots,\gamma_k$ and $\delta_1,\ldots,\delta_l$
  respectively. Then the following are equivalent:
  \begin{enumerate}[(1)]
  \item $\form{c_i}{d_j} = 0$ for all $i, j$

  \item $[T_{\gamma_i},T_{\delta_j}] \in\ker r$ for all $i, j$
  \item $[T_C,T_D]\in\ker r$
  \item $r(T_C)$ and $r(T_D)$ do not generate a free group
  \item $\tr(r(T_CT_D)) = 2g$
\end{enumerate}}

\DeclareMathOperator{\tr}{tr}
\DeclareMathOperator{\Mod}{Mod}
\DeclareMathOperator{\GL}{GL}
\DeclareMathOperator{\Sp}{Sp}

\DeclareMathOperator{\Homeo}{Homeo}
\DeclareMathOperator{\Aut}{Aut}
\DeclareMathOperator{\End}{End}

\newcommand{\iso}{\cong}
\renewcommand{\phi}{\varphi}
\renewcommand{\epsilon}{\varepsilon}

\renewcommand{\S}{S}
\renewcommand{\t}{t}
\newcommand{\fund}{\Gamma}
\newcommand{\MCG}{\Mod_{g,1}}
\newcommand{\arc}{\xi}
\newcommand{\Curve}{\gamma}
\newcommand{\lcurve}{c}

\newcommand{\boundarycurve}{\delta}
\newcommand{\dirac}{\lambda}
\newcommand{\ZH}{\Z H}
\newcommand{\Zfund}{\Z \fund}
\newcommand{\ZN}{\Z N}

\newcommand{\acover}{\widehat{\S}}
\newcommand{\Hcover}{H_1(\acover,\pi^{-1}(*))}
\newcommand{\form}[2]{\langle #1,#2\rangle}
\newcommand{\lift}[1]{\widehat{#1}}
\newcommand{\kerep}{(\ker\epsilon)}
\newcommand{\coloneq}{\mathrel{\mathop:}\mkern-1.2mu=}
\newcommand{\co}{\colon\thinspace}
\newcommand{\mdash}{\hspace{2pt}---\hspace{2pt}}

\begin{document}

\begin{abstract}
  The Magnus representation of the Torelli subgroup of the mapping
  class group of a surface is a homomorphism
  $r\co\mathcal{I}_{g,1}\to\GL_{2g}(\ZH)$. Here $H$ is the first
  homology group of the surface. This representation is not faithful;
  in particular, Suzuki previously described precisely when the
  commutator of two Dehn twists about separating curves is in $\ker
  r$.  Using the trace of the Magnus representation, we apply a new
  method of showing that two endomorphisms generate a free group to
  prove that the images of two positive separating multitwists under
  the Magnus representation either commute or generate a free group,
  and we characterize when each case occurs.

\end{abstract}
\maketitle

\section{Introduction}
Let $\S = S_{g,1}$ be a compact orientable surface of genus $g$ with
one boundary component. The \emph{mapping class group} of $\S$ is
defined to be $\MCG = \pi_0(\Homeo^+(\S))$, where $\Homeo^+(\S)$ is
the group of orientation-preserving self-homeomorphisms of $\S$ fixing
the boundary $\partial\S$ pointwise. In other words, $\MCG$ is the
group of isotopy classes of homeomorphisms of $\S$ relative to the
boundary. The mapping class group naturally acts on $H=H_1(S;\Z)$
preserving the algebraic intersection form; the Torelli group
$\mathcal{I}_{g,1}$ is defined to be the kernel of this
representation, giving the short exact sequence
\[1\to\mathcal{I}_{g,1}\to\MCG\to \Sp_{2g}(\Z)\to 1.\] There are two
basic types of elements in $\mathcal{I}_{g,1}$: a \emph{separating
  twist} is a Dehn twist $T_\Curve$ about a separating simple closed
curve $\Curve$, while a \emph{bounding pair} map is of the form
$T_{\Curve_1}T_{\Curve_2}^{-1}$ where $\Curve_1$ and $\Curve_2$ are
nonseparating simple closed curves which are disjoint and homologous. These suffice
to generate $\mathcal{I}_{g,1}$, by work of Birman \cite{BirmanSiegel}
and Powell~\cite{Powell}.

The Magnus representation was originally
defined algebraically using the Fox free differential calculus. Suzuki \cite{SuzukiGeometric} gave an equivalent
topological definition of the Magnus representation, which we will use
throughout this paper. A homeomorphism of $\S$ can be lifted to a
homeomorphism of the universal abelian cover $\acover$ of $\S$, which
then acts on the relative first homology group
$\Hcover\iso\ZH^{2g}$. Suzuki observed that a classical formula due
to Fox implies that the resulting representation ${r\co \mathcal{I}_{g,1}\to \GL_{2g}(\ZH)}$ coincides with the Magnus representation;
we discuss the details of this construction further in
Section~\ref{section:setup}.

By exploiting the fact that the abelian cover is a surface, we define
a $\ZH$--valued ``higher intersection form'' $\form{\cdot}{\cdot}$
on $\Hcover$ that is preserved by the image of the Magnus
representation. A version of this form was first constructed by
Papakyriakopoulos \cite{Papakyriakopoulos}, Morita verified by
laborious computation that it is preserved by the classical Magnus
representation in \cite{MoritaAbelianQuotients}, and Suzuki
\cite{SuzukiKernel} used the topological definition of the classical
Magnus representation to give a natural proof of Morita's
result. We provide a brief self-contained exposition of the
fundamental properties satisfied by our version of this form.

The higher intersection form can be used to state simple formulas
describing the images of separating twists under the Magnus
representation. We study the trace of the Magnus representation,
defining a class function $\t\co\mathcal{I}_{g,1}\to\ZH$ by 
$\t(f) \coloneq \tr(r(f)) - 2g$. When $f$ is in $\mathcal{K}_{g,1}$,
the subgroup of $\mathcal{I}_{g,1}$ generated by separating twists, it
is possible to compute $t(f)$ explicitly in terms of the higher
intersection form by writing $f$ as a product of separating twists
(see Proposition~\ref{tracecalculation}).

To put our main result in context, we recall some well-known facts
regarding separating twists in $\MCG$. First, two separating twists
$T_\gamma$ and $T_\delta$ commute if and only if $\gamma$ and $\delta$
are disjoint; if $T_\gamma$ and $T_\delta$ do not commute, then they
generate a free group. A \emph{positive multitwist} is a product of
Dehn twists $T_C=T_{\gamma_1}^{n_1}\cdots T_{\gamma_k}^{n_k}$, where
the $\gamma_i$ are pairwise disjoint curves and the $n_i$ are positive
integers. Such a multitwist is $\emph{separating}$ if each $\gamma_i$
is a separating curve. Let $T_C = T_{\gamma_1}^{n_1}\cdots
T_{\gamma_k}^{n_k}$ and $T_D = T_{\delta_1}^{m_1}\cdots
T_{\delta_l}^{m_l}$ be positive separating multitwists. Then $T_C$ and
$T_D$ commute if and only if each $\gamma_i$ is disjoint from each
$\delta_j$; if $T_C$ and $T_D$ do not commute, then they generate a
free group (Hamidi-Tehrani~\cite[Theorem 3.2]{H-T}).

Suzuki \cite[Corollary 4.4]{SuzukiKernel} characterized when the
commutator $[T_{\Curve_1}, T_{\Curve_2}]$ of two Dehn twists around
separating curves $\Curve_1$, $\Curve_2$ is in $\ker r$. We extend
this result in two directions at once by classifying all relations
between the images of two positive separating multitwists under the
Magnus representation. Suzuki's result is the equivalence of (1) and
(2) in the following theorem; the remaining equivalences are
original. In particular, the equivalence of (3) and (4) says that the
images of two positive separating multitwists under the Magnus
representation either commute or generate a free group.

\begin{mainthm}
\tfaetext
\end{mainthm}

\begin{rmk}
  Taking into account the preceding comments, the property
  $\form{c_i}{d_j} = 0$  should be viewed as a weaker version of
  disjointness, which might be called \emph{Magnus-disjointness}. We use
  this notion to define a \emph{Magnus-multitwist} to be a product of
  Dehn twists about Magnus-disjoint curves, and observe that
  Theorem~\ref{thm:main} holds verbatim for positive separating
  Magnus-multitwists.
\end{rmk}

\begin{rmk} The requirement that $T_C$ and $T_D$ be \emph{positive}
  separating multitwists is only necessary to show (5) $\implies$
  (1). The chain of implications \[(1) \implies (2) \implies (3)
  \implies (4) \implies (5)\] holds for arbitrary
  separating multitwists (or Magnus-multitwists) $T_C$ and $T_D$.
\end{rmk}\pagebreak

To prove the main theorem, we use the following result, which may be
of independent interest. It gives conditions under which two
endomorphisms must generate a free group.
\begin{freelem}
\freetext
\end{freelem}

We begin by explaining Suzuki's topological definition of the Magnus
representation in Section~\ref{section:setup}.  In
Section~\ref{section:forms} we define the higher algebraic
intersection form and give a self-contained exposition of its
properties, many of which were previously described by
Papakyriakopoulos \cite{Papakyriakopoulos}, Morita
\cite{MoritaAbelianQuotients}, or Suzuki \cite{SuzukiKernel}. In
Section~\ref{section:separating} we explain the connection between the
Magnus representation and the higher intersection form. In
Section~\ref{section:free} we prove Lemma~\ref{lem:free}. In
Section~\ref{section:main} we combine the preceding results to prove
Theorem~\ref{thm:main}. Finally, in Section~\ref{section:questions},
we outline some natural further questions based on our results.

\medskip\noindent\textbf{Acknowledgements.}
  We thank Masaaki Suzuki for kindly providing us with a preprint of
  \cite{SuzukiGeometric}, which provided the impetus for this
  work. This paper was begun during the REU program at Cornell
  University in 2005 under the supervision of Tara Brendle, funded by
  the National Science Foundation. We thank the REU and Cornell for
  their support and hospitality.  The first author would additionally
  like to thank the Cornell Presidential Research Scholars program,
  whose support in part made possible the completion of this paper.

  We are grateful to Khalid Bou-Rabee, Tara Brendle, Nathan Broaddus,
  Spencer Dowdall, Benson Farb, Asaf Hadari, Ben McReynolds and Justin
  Malestein for reading early versions of this paper, and for their
  valuable comments and advice. We are especially grateful to Matthew
  Day for conversations regarding the contributions of
  Papakyriakopoulos to the subject. We are indebted to two anonymous referees,
  whose careful reading and helpful comments resulted in the
  substantial reorganization of the paper. We also thank Miklos Abert,
  Joan Birman, Ken Brown, Keith Dennis, Allen Hatcher, Martin
  Kassabov, Andrew Putman, and Karen Vogtmann for helpful discussions.

\section{The Magnus representation}\label{section:setup}
The Magnus representation was originally
defined using the Fox free differential calculus
 (see Birman \cite{BirmanBraids}, for instance). In this paper, we
will instead use a topological definition of $r\co
\mathcal{I}_{g,1}\to \GL_{2g}(\ZH)$ that was recently
described by Suzuki in \cite{SuzukiGeometric}, which we now outline.

We fix a basepoint $*\in\partial\S$, and note that the fundamental
group $\fund \coloneq \pi_1(\S,*)$ is free on $2g$ generators. We will
avoid performing any computations with respect to specific generating
sets or bases, but it is occasionally convenient to choose generators
for $\fund$. Let $A_1,\ldots,A_g,B_1,\ldots,B_g$ be a generating set
for $\fund$ such that the product of commutators
$[A_1,B_1]\cdots[A_g,B_g]$ is a positively oriented loop around the boundary
component. Note that if $a_i$, $b_i$ are the homology classes of
$A_i$, $B_i$ respectively, then $a_1,\ldots,a_g$, $b_1,\ldots,b_g$ form a
symplectic basis with respect to the algebraic intersection form on
$H=H_1(\S)$. (Throughout the paper, all homology groups are taken with
coefficients in $\Z$.)

 Let $\pi\co \acover\to\S$ be
the regular covering space corresponding to the commutator subgroup of
$\fund$; this is known as the universal abelian cover of $\S$. The
deck transformations are just $\fund/[\fund,\fund]\iso\Z^{2g}$, and by
identifying this quotient with $H$ the homology groups of $\acover$
are naturally viewed as $\ZH$--modules.

Fix a lift $\widehat{*}\in \pi^{-1}(*)\subset \partial \acover$ of the
basepoint $*$. An element of $\MCG$ can be represented by a
homeomorphism $f\co \S\to\S$ that fixes $\partial\S$ pointwise. Since
$[\fund,\fund]$ is a characteristic subgroup of $\fund$, any
homeomorphism $f$ lifts to a homeomorphism $\lift{f}\co
\acover\to\acover$, and this lift is unique as long we require that
$\lift{f}$ fix $\widehat{*}$.

The action of $\lift{f}$ on the relative
homology $\Hcover$ is obviously $\Z$--linear, but it is not
necessarily $\ZH$--linear\mdash it is twisted by the action of $f$
on $H$. If we restrict $f$ to lie in $\mathcal{I}_{g,1}$, this ensures
the lift $\lift{f}$ will act $\ZH$--linearly on $\Hcover$.  Since
$\Hcover$ is isomorphic to $\ZH^{2g}$ (though this isomorphism is not
canonical, see the discussion below), this action of $\lift{f}$ can be
identified with a representation
\[r\co \mathcal{I}_{g,1}\to \Aut(\Hcover)\iso \GL_{2g}(\ZH),\] which we take as the definition of the Magnus
representation $r$. This allows us to regard $r(f)$ as an automorphism of $\Hcover$, as we will do
throughout the paper. 

We now examine the structure of the module $\Hcover$ in greater
detail. Consider the long exact sequence of homology for the pair
$(\acover,\pi^{-1}(*))$. The only nonzero part of this sequence is
\begin{equation}\label{eq:les}
  0\to H_1(\acover)\to \Hcover
  \overset{\partial}{\longrightarrow} H_0(\pi^{-1}(*))
  \overset{\epsilon}{\longrightarrow} H_0(\acover)\to0,
\end{equation}
where the maps $\epsilon$ and $\partial$ will be described in greater
detail below. As discussed above, the action of $H$ on $\acover$ by
deck transformations makes this an exact sequence of
$\ZH$--modules. Elements of the Torelli group act naturally on the
entire long exact sequence by $\ZH$--module automorphisms.

We describe the $\ZH$--module structure of each term of this
sequence. We begin by noting that $H_0(\acover) = \Z$ because
$\acover$ is connected, while $H_0(\pi^{-1}(*)) \iso \ZH$ because
the connected components of $\pi^{-1}(*)$ correspond to elements of
$H$. Explicitly, ${h\mapsto h(\widehat{*})}$ defines a bijection $H\to
\pi^{-1}(*)$. The map $\epsilon$ in (\ref{eq:les}) clearly corresponds
to the augmentation map $\ZH\to\Z$ that maps $h$ to $1$ for all
$h\in H$.

We can see that $\Hcover\iso \ZH^{2g}$ by lifting a basis for
$\fund$ to a $\ZH$--module basis for $\Hcover$, as follows. For each
$i=1,\ldots,g$, define $\alpha_i\in\Hcover$ to be the homology class
of the unique lift of the loop $A_i$ starting at $\widehat{*}$, and
similarly let $\beta_i$ be the lift of $B_i$. Each of these arcs must
have its endpoints in $\pi^{-1}(*)$, so each describes an element of
$\Hcover$. The deck transformation $a_i$ translates the tail of the
arc $\alpha_i$ (that is, $\widehat{*}$) to its head, and similarly
$b_i$ translates $\widehat{*}$ to the head of $\beta_i$. Note that $\acover$ deformation retracts onto the 1--dimensional spine whose edges are the translates of $\alpha_1,\ldots,\alpha_g, \beta_1,\ldots,\beta_g$. Thus $\Hcover$ is a free $\ZH$--module of rank $2g$,
with basis $\alpha_1,\ldots,\alpha_g, \beta_1,\ldots,\beta_g$. It is
easily verified that the map $\partial$ in (\ref{eq:les}) is given by
$\partial \alpha_i = a_i-1$ and $\partial \beta_i = b_i-1$. The
isomorphism $\Hcover \iso \ZH^{2g}$ clearly is non-canonical: the
choice of basis corresponds to the choice of generators of $\fund$.\pagebreak

Finally, we have $\pi_1(\acover)= [\fund,\fund]$, so $H_1(\acover)$ is
the abelianization $[\fund,\fund]^{\text{ab}}$. The action of $\fund$
on $[\fund,\fund]$ by conjugation descends to a action of $H$ on
$[\fund,\fund]$ by outer automorphisms. This projects to a
well-defined action of $H$ on $[\fund,\fund]^{\text{ab}}$, and the
resulting $\ZH$--module structure agrees with that of
$H_1(\acover)$. The exact sequence above is thus isomorphic as an
exact sequence of $\ZH$--modules to
\begin{equation*}
  0\to [\fund,\fund]^{\text{ab}} \to \ZH^{2g}
  \overset{\partial}{\longrightarrow} \ZH
  \overset{\epsilon}{\longrightarrow} \Z \to0.
\end{equation*}

Observe that the pre-image under $\pi$ of any separating curve
$\Curve$ in $S$ is a disjoint union of curves in $\acover$,
each a lift of $\Curve$, whose union separates $\acover$;
 by abuse of notation, we also call the
homology classes of these curves the \emph{lifts} of $\Curve$.
Distinct lifts of $\Curve$ differ by a deck transformation. Note that
any lift of a separating curve, considered as an element of $\Hcover$,
lies in the image of $H_1(\acover)$, which is also the kernel of the
boundary map $\partial$.

\section{Higher intersection forms}
\label{section:forms}

Inspired by Reidemeister, Papakyriakopoulos \cite{Papakyriakopoulos}
defined a biderivation on $\Zfund$ that encodes the intersection
theory of the regular covering spaces of a surface. Given a normal
subgroup $G$ of $\Gamma$ with quotient $N = \Gamma/G$,
Papakyriakopoulos' biderivation yields a $\ZN$-valued pairing on the
first homology of the corresponding regular covering space. When $G =
\Gamma$, this is simply the algebraic intersection form on the base
surface. In general, the pairing can be easily expressed in terms of
the usual  $\Z$--valued algebraic intersection form on the covering
space; see Turaev \cite{Turaev}, or Perron \cite{Perron} for a modern treatment. The case when
$G = [\Gamma,\Gamma]$, so that $N = H$ and the covering space is the
universal abelian cover $\acover$, is of particular interest. Hempel
\cite{Hempel} made frequent use of this pairing on $H_1(\acover)$,
interpreting it as the ``Reidemeister pairing'' defined by
Reidemeister in \cite{Reidemeister}. Suzuki \cite{SuzukiKernel} was
the first to use this pairing (which we will call the \emph{higher
  intersection form}) to study the Magnus representation. He viewed it
as a pairing of separating curves on the base surface $S$, and he
proved that the commutator of two separating twists is in $\ker r$ if
and only if the higher intersection number of the two curves involved
is zero.

The higher intersection form will be intrinsic to our study of the
Magnus representation, so we will find it convenient to use a more
general version of this pairing. In the previous section, we saw that
the Magnus representation can be viewed as the action of the Torelli
group on the $\ZH$-module $\Hcover$, and so we want to define the
higher intersection form on this module rather than on the submodule
$H_1(\acover)$. This is not difficult, and indeed can be done just
using Papakyriakopoulos' biderivation, but we prefer a more geometric
approach.

In this section, we provide a self-contained exposition of the
definition and fundamental properties of the higher intersection form
on $\Hcover$. Many of the results in this section have been previously
proven by Papakyriakopoulos, Suzuki, and others. Other results are
well-known but have been used without proof; we include proofs of such
results for completeness.\pagebreak

We first need to define a $\Z$--valued intersection number of two
elements of  $\Hcover$. Any element of $\Hcover$ can be realized as a
linear combination of immersed closed curves in $\acover$ and immersed
arcs in $\acover$ with endpoints in $\pi^{-1}(*)$. Any pair of curves,
or a curve and an arc, can be realized so that they only intersect
transversely, and then the orientation of $\acover$ gives a natural
algebraic intersection number. This cannot be extended to the case of a pair of arcs,
however; $\alpha_1$ and $\beta_1$, for example, can be represented by
two arcs in $\acover$ that share one endpoint, so their intersection
number is not well-defined. In order to define their algebraic
intersection number, we need to move the basepoint of one arc
slightly; there are two different ways of doing so, and we will keep
track of the resulting differences in the algebraic intersection
number. This will give us two $\Z$--valued bilinear forms on
$\Hcover$.

We now formalize the above discussion. The orientation on $\S$
determines an orientation on $\partial\S$. Take $*'\ne *$ to be a
second basepoint in $\partial\S$; there are two arcs from $*$ to $*'$
contained in $\partial\S$. Call these two arcs $\Curve_+$ and
$\Curve_-$, where $\Curve_+$ is the arc that is positively oriented.
For $\sigma \in \{+,-\}$ let $\arc_\sigma$ be the lift of
$\Curve_\sigma$ to $\acover$ based at $\widehat{*}$.  Then we have
isomorphisms $\phi_\sigma\co \Hcover\to H_1(\acover, \pi^{-1}(*'))$
defined by $\phi_\sigma(x)=x+\partial{x}\cdot\arc_\sigma$. Here $\partial{x}\in \ZH$ acts on the arc $\arc_\sigma$ by translation.
 Thus
$\phi_\sigma$ slides the basepoint along the boundary from
$\pi^{-1}(*)$ to $\pi^{-1}(*')$: $\phi_+$ in the positive direction,
$\phi_-$ in the negative.

The orientation on $\S$ induces an orientation of $\acover$. This
orientation determines a  $\Z$--bilinear algebraic intersection form on
$\Hcover\times H_1(\acover, \pi^{-1}(*'))$, since  representative arcs
have distinct endpoints and thus can always be made transverse.  We denote
this form by $(\cdot,\cdot)$ and define two bilinear forms
$(\cdot,\cdot)_\sigma$ on  $\Hcover\times\Hcover$ by
$(c,d)_\sigma = (c,\phi_\sigma(d))$ for $\sigma\in\{+,-\}$.

We use the action of deck transformations on these forms to define new forms
$\form{\cdot}{\cdot}_\sigma$ below. The formula
(\ref{intersectionform}) below appears in the work of
Papakyriakopoulos; his results \cite[Theorem 10.13]{Papakyriakopoulos}
imply that the intersection form $\form{\cdot}{\cdot}_+$ is equivalent
to the biderivation mentioned in the introduction to this section (and
the other intersection form corresponds to a slightly altered
biderivation).
\begin{definition}
  The two higher intersection forms $\form{\cdot}{\cdot}_\sigma$ are
  $\Z$--bilinear functions
  \[\form{\cdot}{\cdot}_\sigma\co\Hcover\times\Hcover\to\ZH\] for
  $\sigma\in\{+,-\}$, defined by
\begin{equation}\label{intersectionform}
  \form{c}{d}_\sigma = \sum_{h\in H}{(c,hd)_\sigma h}.
\end{equation}
\end{definition}
Note that this sum is finite because $c$ and $d$ are compactly
supported. Morita's explicit calculations in the proof of
\cite[Theorem 5.3]{MoritaAbelianQuotients} imply that the Magnus
representation preserves the forms $\form{\cdot}{\cdot}_\sigma$. In
fact, it is easy to see directly that these forms are preserved by a
general class of topologically defined automorphisms of $\Hcover$
which includes those defining the Magnus representation. Suzuki used
this approach in \cite{SuzukiGeometric} to give the first natural
proof of Morita's result.\pagebreak

\begin{lemma}\label{lem:preserve}
  Suppose that $f$ is an orientation-preserving self-homeomorphism of
  the pair $(\acover,\pi^{-1}(*))$ that commutes with all deck
  transformations $h\in H$. Then the action of $f$ on $\Hcover$
  preserves both forms $\form{\cdot}{\cdot}_\sigma$. In particular,
  $\form{\cdot}{\cdot}_\sigma$ is preserved by the action of $H$ and
  by the image of the Magnus representation $r$.
\end{lemma}
\begin{proof}
  Using the fact that the intersection forms $(\cdot,\cdot)_\sigma$
  are preserved by such a homeomorphism $f$, we have that
  \begin{align*}
    \form{f(c)}{f(d)}_\sigma &= \sum_{h\in H}{(f(c),hf(d))_\sigma h} \\
    &= \sum_{h\in H}{(f(c),f(hd))_\sigma h} \\
    &= \sum_{h\in H}{(c,hd)_\sigma h} = \form{c}{d}_\sigma,
  \end{align*}
  as desired.
\end{proof}

The following lemma corresponds to the fact that Papakyriakopoulos'
pairing mentioned above is a biderivation.  Recall that the
conjugation involution 
$\overline{\hspace{.2em}\vphantom{h}\cdot\hspace{.15em}}\co
\ZH\to\ZH$ is defined for $h\in H$ by $\overline{h} = h^{-1}$ and
on $\ZH$ by extending $\Z$--linearly.

\begin{lemma}\label{lem:sesqui}
  Both forms $\form{\cdot}{\cdot}_\sigma$ are $\ZH$--sesquilinear,
  meaning that
  \begin{equation*}
    \form{gc}{hd}_\sigma = g\overline{h}\form{c}{d}_\sigma
  \end{equation*}
  for any $g,h\in\ZH$ and $c,d\in\Hcover$.
\end{lemma}
\begin{proof}
  By $\Z$--bilinearity, it is sufficient to prove this result when
  $g,h\in H$. From the definition of $\form{\cdot}{\cdot}_\sigma$, we
  have that
  \begin{equation*}
    \form{gc}{hd}_\sigma = \form{c}{g^{-1}hd}_\sigma
    = \sum_{j\in H} (c,jg^{-1}hd)_\sigma j = \sum_{j\in H} (c,jd)_\sigma gh^{-1}j
    = gh^{-1}\form{c}{d}_\sigma
  \end{equation*}
  since $H$ is abelian. The first equality holds because
  $\form{\cdot}{\cdot}_\sigma$ is preserved by the action of $g$
  (Lemma~\ref{lem:preserve}).
\end{proof}

Let $\pi_{\ast}\co \Hcover\to H$ denote the map on first homology
induced by the covering map $\pi\co \acover\to\S$. Recall that
$\epsilon\co\ZH\to\Z$ is the augmentation map. Our justification for
calling these forms $\form{\cdot}{\cdot}_\sigma$ higher intersection
forms is the following lifting property. This property is implicit in
the work of Papakyriakopoulos (for example, it can be derived from
\cite[Formula 8.5]{Papakyriakopoulos}). Here we denote the $\Z$--valued
algebraic intersection form on $H$ by $(\cdot,\cdot)$.

\begin{lemma}\label{lifting}
  For any $c,d\in\Hcover$ and any $\sigma\in\{+,-\}$ we have
  \[\epsilon(\form{c}{d}_\sigma) = (\pi_*(c),\pi_*(d)).\]
\end{lemma}
\begin{proof}
  We want to show that $\sum_{h\in H}(c,hd)_\sigma =
  (\pi_*(c),\pi_*(d))$ for $c,d\in\Hcover$. Since
  $\pi_*(\phi_\sigma(d)) = \pi_*(d)$ in $H_1(S)$, this is equivalent
  to the claim that  $\sum_{h\in H}(c,hd) = (\pi_*(c),\pi_*(d))$ for
  $c\in\Hcover$ and $d\in H_1(\acover,\pi^{-1}(*'))$.

  Since the projected $1$--cycles $\pi_*(c)$ and $\pi_*(d)$ have
  distinct basepoints, we can realize $c$ and $d$ such that $\pi_*(c)$
  and $\pi_*(d)$ intersect transversely in $S$. Then it is easy to see
  that each intersection of $\pi_*(c)$ and $\pi_*(d)$ corresponds to
  an intersection of $c$ with some translate of $d$: if $\pi_*(c)$ and
  $\pi_*(d)$ intersect at some point $x\in S$, then exactly one
  element of $\pi^{-1}(x)$ is in the image of $c$, and exactly one
  translate of $d$ passes through that element. This yields the
  desired result.
\end{proof}

This lifting property of the higher intersection forms allows us to
deduce the nondegeneracy of $\form{\cdot}{\cdot}_\sigma$ from the
nondegeneracy of the symplectic intersection form on $H$. In fact, we
obtain the following stronger result. Note that $\kerep^n$ is an ideal of $\ZH$, and $\kerep^n\Hcover$ is the submodule of $\Hcover$ spanned by terms of the form $r\cdot c$ for $r\in \kerep^n$ and $c\in\Hcover$.

\begin{proposition}\label{nondegenerate}
  For any $n\ge 0$, $\sigma\in\{+,-\}$, and $y\in\Hcover$, we have
  \[y\in\kerep^n\Hcover\quad\iff\quad\form{x}{y}_\sigma\in\kerep^n\text{
    for all }x\in\Hcover.\]
\end{proposition}
\begin{proof}
  If $y\in\kerep^n\Hcover$, it is immediate that
  $\form{x}{y}_\sigma\in\kerep^n$ for any $x$. We prove the other
  implication by induction on $n$; this is trivial if $n = 0$, so
  assume that $n>0$ and that the proposition holds for all smaller values of
  $n$. Suppose for contradiction that $y\notin\kerep^n\Hcover$ but
  $\form{x}{y}_\sigma\in\kerep^n$ for all ${x\in\Hcover}$. Let
  $\{s_1,\ldots,s_{2g}\}$ be the standard basis
  $\{\alpha_1,\ldots,\alpha_g,\beta_1,\ldots,\beta_g\}$ for
  $\Hcover\iso\ZH^{2g}$, taken in any order. By the inductive
  hypothesis, ${y\in\kerep^{n-1}\Hcover}$, so we can write $y =
  \sum_{j=1}^{2g}h_js_j$, where each\linebreak ${h_j\in\kerep^{n-1}}$ but
  $h_{j_0}\notin\kerep^n$ for some $j_0$.

  Now, using Lemma~\ref{lifting}, the nondegeneracy of the
  intersection form on $H$ implies that we can find $x\in\Hcover$ such
  that $\form{x}{s_j}_\sigma\in\kerep\iff j\ne j_0$. Then \[
    \form{x}{y}_\sigma =
    \sum_{j=1}^{2g}\overline{h}_j\form{x}{s_j}_\sigma
    \in\kerep^n\iff\overline{h}_{j_0}\form{x}{s_{j_0}}_\sigma
    \in\kerep^n.\]We have\[
    \overline{h}_{j_0}\form{x}{s_{j_0}}_\sigma\in\kerep^n\iff
    \overline{h}_{j_0}\epsilon(\form{x}{s_{j_0}}_\sigma) \in\kerep^n
    \] because the difference of the two expressions is
  \[\overline{h}_{j_0}\big(\form{x}{s_{j_0}}_\sigma
  -\epsilon(\form{x}{s_{j_0}}_\sigma) \big)\in \kerep^{n-1}\kerep.\]
  But $\overline{h}_{j_0}\epsilon(\form{x}{s_{j_0}}_\sigma)
  \not\in\kerep^n$ because $\overline{h}_{j_0}\notin\kerep^n$ and
  $\form{x}{s_{j_0}}\not\in\kerep$, so we conclude that
  $\form{x}{y}_\sigma\not\in \kerep^n$. This contradiction completes
  the induction.
\end{proof}

Since $\bigcap_{n\ge 0}\kerep^n\Hcover = 0$, we have the following
corollary:
\begin{corollary}\label{nondeg}
  Both forms $\form{\cdot}{\cdot}_\sigma$ are nondegenerate; that is,
  for each ${\sigma\in \{+,-\}}$ and every ${x\in \Hcover\backslash\{0\}}$, there
  exists $y\in \Hcover$ such that $\form{x}{y}_\sigma\ne 0$.
\end{corollary}

The following proposition tells us the difference between the two
higher intersection forms.

\begin{proposition}\label{formdifference}
  Suppose $c,d\in \Hcover$. Then $\form{c}{d}_+-\form{c}{d}_-
  =\partial c \overline{\partial d}\in\ZH$.
\end{proposition}
\begin{rmk} When either $c$ or $d$ is in $\ker\partial$, Proposition~\ref{formdifference} shows
  that $\form{c}{d}_+ =\form{c}{d}_-$, so in this case we need only
  write $\form{c}{d}$. In particular, this is the case whenever $c$ or
  $d$ is the lift of a separating curve.
\end{rmk}\pagebreak
\begin{proof} Let $\mathbf{\boundarycurve}=\arc_+-\arc_-$ be a loop around the
  boundary (the lift of $\partial\S$ starting at $\widehat{\ast}$). Note that $(c,\delta)$ measures how many times $c$ crosses $\delta$, which we can identify with the coefficient of 1 in $\partial c$. Similarly, for a translate $g\delta$, if we write $\partial c=\sum_{h\in H}\alpha_h h$ we have $(c,g\delta)=\alpha_g$. This implies that \[\form{c}{\delta}=\sum_{h\in H}(c,h\delta)h=\sum_{h\in H}\alpha_h h=\partial c.\] The definition of $\delta$ shows that
  \begin{align*}
    (x,y)_+-(x,y)_-&=\big(x,\phi_+(y)\big)-\big(x,\phi_-(y)\big)\\
    &=\big(x,(y+\partial y\cdot\arc_+)-(y+\partial y\cdot\arc_-)\big)
    \\&=\big(x,\partial y(\arc_+-\arc_-)\big) =\big(x,(\partial
    y)\boundarycurve\big).
  \end{align*} Applying this identity to each term of $\form{c}{d}_+-\form{c}{d}_-$, we obtain:
\begin{align*}\form{c}{d}_+-\form{c}{d}_-&=\sum_{h\in H}{(c,hd)_+ h}-\sum_{h\in
      H}{(c,hd)_- h}\\&=\sum_{h\in H}(c,h(\partial d)\delta)h\\&=\form{c}{(\partial d)\delta}\\&=\overline{\partial d}\form{c}{\delta}=\overline{\partial d}\partial c\end{align*} as desired.
\end{proof}

We also have the following ``antisymmetry'' property of
$\form{\cdot}{\cdot}_\sigma$:

\begin{lemma}\label{antisym}
  Suppose $\sigma \in \{+,-\}$ and $c,d\in\Hcover$. Then
  ${\form{d}{c}_\sigma = -\overline{\form{c}{d}_{-\sigma}}}$.
\end{lemma}
\begin{proof}
  Note that $(f,e)_\sigma = -(e,f)_{-\sigma}$ for any $e,f\in\Hcover$. Thus
  \begin{align*}
    \form{d}{c}_\sigma &= \sum_{h\in H}{(d,hc)_\sigma h} \\ &=
    -\sum_{h\in H}{(hc,d)_{-\sigma} h} \\ &= -\sum_{h\in
      H}{(c,h^{-1}d)_{-\sigma} h} \\ &= -\sum_{h\in
      H}{(c,hd)_{-\sigma} h^{-1}} =
    -\overline{\form{c}{d}_{-\sigma}}.\qedhere
  \end{align*}
\end{proof}

  \begin{rmk} When $c$ or $d$ is in $\ker\partial$, this shows that
    $\form{d}{c} = -\overline{\form{c}{d}}$, so 
    ${\form{c}{d}=0}\iff{\form{d}{c} = 0}$.  Note that this lemma does
    not imply that $\form{c}{c} = 0$ if $c\in\ker\partial$, but only
    that ${\form{c}{c}=-\overline{\form{c}{c}}}$. However, we do have
    the weaker statement that $\form{c}{c} = 0$ if $c$ is a lift of a
    separating curve in the base surface $\S$, as the various lifts
    $\left\{hc\left|\,h\in H\right.\right\}$ are disjoint in this
    case.
\end{rmk}\pagebreak

\section{The restriction of $r$ to $\mathcal{K}_{g,1}$}\label{section:separating}
The higher intersection form of the previous section is of great use
in describing the image of a product of separating twists under the
Magnus representation $r$. In this section, we develop this connection
to further understand the restriction of the Magnus representation to
the Johnson kernel $\mathcal{K}_{g,1}$, the subgroup of the Torelli
group generated by separating twists.

The following fundamental result illustrates the relationship between
the higher intersection form and the Magnus representation. This
formula is claimed without proof in Suzuki~\cite{SuzukiKernel} (though
it is stated incorrectly there).
\begin{proposition}\label{prop:twistformula}
  Let $\Curve$ be a separating curve in $\S$ and let
  $\lcurve\in\Hcover$ be any lift of $\Curve$. Let $T_\Curve$ denote
  the Dehn twist around $\Curve$. Then the action of $\lift{T}_\Curve$ on
  $x\in\Hcover$ is given by
  \begin{equation*}
    r(T_\Curve)(x)=x+\form{x}{\lcurve}\lcurve.
  \end{equation*}
\end{proposition}

\begin{rmk}
  Recall that the lift $\lcurve$ lies in $\ker\partial$
  exactly when $\Curve$ is a separating curve in $\S$, so $\form{x}{\lcurve} =
  \form{x}{\lcurve}_+ = \form{x}{\lcurve}_-$ in the above formula.
  The result is analogous to the formula for the action of a Dehn
  twist on $H$ given by ${T_\Curve(h)=h+(h,[\Curve])[\Curve]},$ where
  $(\cdot,\cdot)$ is the algebraic intersection form on $H$.
 This may be taken as evidence that this definition of
  $\form{\cdot}{\cdot}$ is the correct one.
\end{rmk}
\begin{proof}
  The lifted homeomorphism $\lift{T}_\Curve$ can be thought of as
  simultaneously twisting about each lift of $\Curve$, since
  these lifts are nonintersecting closed curves in $\acover$. For each
  intersection of $x$ with a lift $\lift{\Curve}$ of $\Curve$, we add
  or subtract $\lift{\Curve}\in\Hcover$, depending on the orientation
  of the intersection. These lifts of $\Curve$ are simply the curves
  $hc$ for $h\in H$. Thus
  \begin{equation*}
    r(T_\Curve)(x)=x+\sum_{\lift{\Curve} \text{ lifts }\Curve}
    (x,\lift{\Curve})\lift{\Curve}=x+\sum_{h\in
      H}(x,h\lcurve)h\lcurve.
  \end{equation*}
  By equation (\ref{intersectionform}), this is just
  $x+\form{x}{\lcurve}\lcurve$.
\end{proof}

Recall that a separating multitwist is a product of Dehn twists $T_C =
T_{\gamma_1}^{n_1}\cdots T_{\gamma_k}^{n_k}$ such that each $\gamma_i$
is a separating curve, pairwise disjoint. We
generalize this concept by defining a
\emph{separating Magnus-multitwist} to be a product of Dehn twists ${T_C =
T_{\gamma_1}^{n_1}\cdots T_{\gamma_k}^{n_k}}$ such that each $\gamma_i$
is a separating curve with lift $c_i$ satisfying 
$\form{c_i}{c_{i'}}=0$ for all $i,i'$.  (For example, if $\gamma_i$
and $\gamma_{i'}$ have geometric intersection number $2$, then we have
$\form{c_i}{c_{i'}}=0$.) Note that a separating multitwist is also a
separating Magnus-multitwist: if the $\gamma_i$ are disjoint, it
follows that their lifts $c_i$ are disjoint, and so
$\form{c_i}{c_{i'}}=0$ for all $i,i'$. Then
Proposition~\ref{prop:twistformula} yields the following corollary.
\begin{corollary}\label{cor:multitwistformula}
  Let $T_C=T_{\gamma_1}^{n_1}\cdots T_{\gamma_k}^{n_k}$ be a
  separating Magnus-multitwist. Then the action of $T_C$ on
  $x\in\Hcover$ is given by
  \begin{equation*}
    r(T_C)(x)=x+\sum_{i = 1}^{k}n_i\form{x}{c_i}c_i.
  \end{equation*} 
\end{corollary}\pagebreak

In the remainder of this section, we analyze the trace of the Magnus
representation $r$. Let $\t\co \mathcal{I}_{g,1}\to\ZH$ be defined by
$\t(f) \coloneq \tr (r(f))-2g$; this is a class function on the
Torelli group.  The normalization is chosen such that $\t(1) = 0$.  If
$T_C$ is a separating multitwist, it follows from
Corollary~\ref{cor:multitwistformula} that $r(T_C)-1$ is nilpotent and
thus has trace 0. As a consequence, for any separating multitwist
$T_C$, we have ${\t(T_C)=\tr(r(T_C)-1)=0}$.

In general, we can express the value
of $t$ on $\mathcal{K}_{g,1}$ in terms of the higher intersection form:
\begin{proposition}\label{tracecalculation}
  Let $\gamma_1,\ldots,\gamma_k$ be separating curves in $\S$ with
  lifts $c_1,\ldots,c_k$ in  $\Hcover$, and let
  $n_1,\ldots,n_k$ be integers. Then
  \begin{equation*}
    \t(T_{\Curve_1}^{n_1}T_{\Curve_2}^{n_2}\cdots T_{\Curve_k}^{n_k})
    = \sum_{\substack{1\le i_1 < i_2 < \cdots < i_m \le k \\ m\ge 2}}
    n_{i_1}n_{i_2}\cdots n_{i_m}\form{\lcurve_{i_1}}{\lcurve_{i_m}}
    \form{\lcurve_{i_m}}{\lcurve_{i_{m-1}}}\cdots\form{\lcurve_{i_2}}{\lcurve_{i_1}}.
  \end{equation*}
\end{proposition}
\begin{proof}
  Use Proposition~\ref{prop:twistformula} to expand
  $r(T_{\Curve_1}^{n_1}T_{\Curve_2}^{n_2}\cdots
  T_{\Curve_k}^{n_k})(x)$ as a sum of $2^k$ terms. Since taking the
  trace is linear, we can then write
  $\t(T_{\Curve_1}^{n_1}T_{\Curve_2}^{n_2}\cdots T_{\Curve_k}^{n_k})$
  as
  \[\sum_{\substack{1\le i_1 < i_2 < \cdots < i_m \le k \\ m\ge 1}}
  n_{i_1}n_{i_2}\cdots n_{i_m}\tr\!\big[x\mapsto \form{x}{\lcurve_{i_m}}
  \form{\lcurve_{i_m}}{\lcurve_{i_{m-1}}}\cdots
  \form{\lcurve_{i_2}}{\lcurve_{i_1}}\lcurve_{i_1}\big],\] where the term
  $\tr[x\mapsto x]$ has canceled out the $-2g$.  Taking the trace of a
  rank $1$ operator is easy: $\tr[x\mapsto \lambda(x)v] = \lambda(v)$ for any
  linear functional $\lambda:\ZH^{2g}\to \ZH$ and element
  $v\in\ZH^{2g}$.
  Using this fact and observing that the terms with $m = 1$ vanish
  because $\form{c}{c} = 0$ for any lifting curve $c$, we obtain the
  desired result.
\end{proof}

\begin{rmk}
  The special case $k = 2$, $n_1 = n_2 = 1$ of
  Proposition~\ref{tracecalculation} was previously obtained by Suzuki
  \cite[Theorem 4.3]{SuzukiKernel} via different methods.
\end{rmk}

The properties of the Fox derivatives imply that for any
$f\in\mathcal{K}_{g,1}$, the entries of $r(f)-1$ lie in $\kerep^2$ and
so $t(f)\in\kerep^2$. For such $f$, Morita's trace functional is the
homomorphism sending $f$ to $[t(f)]\in \kerep^2/\kerep^3$. The
vanishing of this trace functional, proved by Morita in
\cite{MoritaAbelianQuotients}, is equivalent to the statement that in
fact $t(f)\in \kerep^3$ for all $f\in \mathcal{K}_{g,1}$. Using
Proposition~\ref{tracecalculation}, we can extend this result as
follows.
\begin{corollary}\label{cor:kerep}
  For any $f\in\mathcal{K}_{g,1}$, $\t(f)\in\kerep^4$.
\end{corollary}
\begin{proof}
  It is easily verified that $\ker\partial$ is contained in the
  submodule  $\kerep\Hcover$. Thus by Lemma~\ref{lem:sesqui},
  $\form{c}{d}\in\kerep^2$ for any $c,d\in \ker\partial$. The desired
  result then follows, since each term in the formula given by
  Proposition~\ref{tracecalculation} contains at least two such
  factors.
\end{proof}

\pagebreak
\section{Showing two endomorphisms generate a free group}
\label{section:free}

The following lemma gives a new method of showing that two
endomorphisms generate a free group.

\begin{lemma}\label{lem:free}
  \freetext
\end{lemma}
\begin{proof}
  Suppose for contradiction that there is some nontrivial relation
  between $1+A$ and $1+B$. Note that $(1+A)^n=1+nA$ and
  $(1+B)^m=1+mB$. Since $\tr(AB)$ is transcendental, in particular $A$
  and $B$ must be nonzero, and thus $(1+A)^n\neq 1$ and $(1+B)^n\neq
  1$ for nonzero $n\in\Z$.  Any nontrivial relation thus involves both
  $1+A$ and $1+B$, and so must be conjugate to one of the form
  \begin{equation*}
    (1+B)^{m_1}(1+A)^{n_1}\cdots(1+B)^{m_k}(1+A)^{n_k} = 1
  \end{equation*}
  for $k > 0$ and nonzero $m_1,n_1,\ldots,m_k,n_k\in\Z$. Multiplying
  on the left by $A$ and on the right by $B$, we have
  \begin{equation*}
    A(1+m_1B)(1+n_1A)\cdots(1+m_kB)(1+n_kA)B-AB=0.
  \end{equation*}
  Since $A^2 = B^2 = 0$, when we expand this expression every term
  vanishes except those of the form $nAB\cdots AB$. Thus we can write
  this expression as\linebreak $P(AB)=0$, where $P(X)$ is a polynomial with
  integer coefficients and leading term $m_1n_1\cdots
  m_kn_kX^{k+1}$. Then all the characteristic values of $AB$ are roots
  of the polynomial $P$, so their sum $\tr(AB)$ is algebraic over
  $\Q$. This contradicts the hypothesis, so no such relation between
  $1+A$ and $1+B$ can exist.  \end{proof}

\section{Proof of the main theorem}
\label{section:main}
In this section, we prove our main theorem, which classifies relations
between the images of two separating multitwists under the Magnus
representation. As mentioned in the introduction, the equivalence of
(1) and (2) was proved by Suzuki \cite{SuzukiKernel}; the remaining
equivalences are original.

\begin{theorem}\label{thm:main}
  \tfaetext
\end{theorem}
\begin{proof}
  (1) $\implies$ (2): If $\form{c_i}{d_j}=0$, then
  $[T_{\gamma_i},T_{\delta_j}]$ is a separating Magnus-multitwist, so
  Corollary~\ref{cor:multitwistformula} gives
  \[r([T_{\gamma_i},T_{\delta_j}])(x)
  =x+\form{x}{c_i}c_i+\form{x}{d_j}d_j-\form{x}{c_i}c_i-\form{x}{d_j}d_j=x\]
  and thus $[T_{\gamma_i},T_{\delta_j}]\in \ker r$, as desired.\pagebreak

  (2) $\implies$ (3): $r(T_C)$ and $r(T_D)$ are products of the pairwise
  commuting elements $r(T_{\gamma_i})$ and $r(T_{\delta_j})$, so they also commute.\\

  (3) $\implies$ (4): $[r(T_C),r(T_D)]=1$ is a nontrivial relation between
  $r(T_C)$ and $r(T_D)$.\\

  (4) $\implies$ (5): Let $K$ be the field of fractions of the
  integral domain $\ZH$. Any $\ZH$-linear endomorphism of
  $\Hcover$ extends $K$-linearly to an endomorphism of the $K$-vector
  space $V \coloneq \Hcover\bigotimes_{\ZH}K$, and in this way we will
  treat elements of the image of the Magnus representation as
  endomorphisms of $V$. We apply Lemma~\ref{lem:free} to the
  endomorphisms $A=r(T_C)-1$ and $B=r(T_D)-1$ of $V$; it follows from
  Corollary~\ref{cor:multitwistformula} that $A^2 = B^2 =
  0$. Lemma~\ref{lem:free} then implies that if $r(T_C)$ and $r(T_D)$
  do not generate a free group, then $\tr(AB)\in \ZH\subseteq K$ is
  algebraic (over $\Q$). The algebraic elements of $\ZH$ are
  precisely those in $\Z$, so this means that $\tr(AB)\in \Z$. Since
  $A$ and $B$ are nilpotent, $\tr(A)=\tr(B)=0$, so
  $\tr(r(T_C))=\tr(r(T_D))=2g$. Now
  expanding \[AB=[r(T_C)-1][r(T_D)-1]=r(T_CT_D)-r(T_C)-r(T_D)+1\] and
  taking traces, we get that
  $\tr(AB)=\tr\big(r(T_CT_D)\big)-2g=\t(T_CT_D)$. By
  Corollary~\ref{cor:kerep}, $\t(T_CT_D)\in\kerep^4\le\ker\epsilon$,
  so $\t(T_CT_D)=\tr(AB)\in\Z$ implies that $t(T_CT_D) = 0$, and thus
  $\tr(r(T_CT_D)) = 2g$.\\

  (5) $\implies$ (1): Using Proposition~\ref{tracecalculation} and the
  fact that $\form{c_i}{c_{i'}} = \form{d_j}{d_{j'}} = 0$ for any
  $i,i'$ or $j,j'$, we calculate that
  \begin{equation*}
    0 = \tr(r(T_CT_D))-2g = t(T_CT_D) = \sum_{\substack{1\le i \le k \\ 1 \le j \le l}}
    n_im_j\form{c_i}{d_j}\form{d_j}{c_i}.
  \end{equation*}
  By Lemma~\ref{antisym}, we can rewrite this as
  \begin{equation*}
    \sum_{\substack{1\le i \le k \\ 1 \le j \le l}}
    n_im_j\form{c_i}{d_j}\overline{\form{c_i}{d_j}} = 0.
  \end{equation*}
  If we define for each $g\in H$ a linear functional $\dirac_g\co \ZH\to \Z$ by
  $\dirac_g\!\big(\sum_{h\in H}\alpha_h h\big)=\alpha_g$, note in general
  that \[\dirac_1(x\overline{x})=\sum_{h\in H}(\dirac_h x)^2\] is
  nonnegative, and that $\dirac_1(x\overline{x})=0$ implies
  $x=0$. Applying $\dirac_1$ to the equation above yields
  \begin{equation*}
    \sum_{\substack{1\le i \le k \\ 1 \le j \le l}}
    n_im_j \dirac_1\!\big(\form{c_i}{d_j}\overline{\form{c_i}{d_j}}\big) = 0.
  \end{equation*}
  Since the $n_i$ and $m_j$ are all positive, this implies that
  $\dirac_1\!\big(\form{c_i}{d_j}\overline{\form{c_i}{d_j}}\big) = 0$ for
  all $i$ and $j$; thus $\form{c_i}{d_j}=0$ for all $i$ and $j$, as
  desired.
\end{proof}\pagebreak

\begin{rmk}
  In the proof of Theorem~\ref{thm:main}, we use only that
  $\form{c_i}{c_{i'}}=\form{d_j}{d_{j'}}=0$. Recalling that we call
  $T_C=T_{\gamma_1}^{n_1}\cdots T_{\gamma_k}^{n_k}$ a
  \emph{separating Magnus-multitwist} if each $\gamma_i$ is a
  separating curve and the lifts $c_i$ of $\gamma_i$ satisfy
  $\form{c_i}{c_{i'}}=0$ for any $i,i'$, we have the following corollary.
\end{rmk}
\begin{corollary}
  Let $T_C$ and $T_D$ be positive separating Magnus-multitwists. Then
   Theorem~\ref{thm:main} holds verbatim.
\end{corollary}

\begin{rmk}
  As mentioned in the introduction, the only case in the proof
  of Theorem~\ref{thm:main} that requires $T_C$ and $T_D$ to
  be \emph{positive} separating multitwists was {(5) $\implies$
    (1)}. The other implications hold for arbitrary separating
  multitwists (or Magnus-multitwists) $T_C$ and $T_D$.
\end{rmk}

\section{Further questions}
\label{section:questions}

There are several natural ways in which one might try to generalize
the classification described in this paper; we describe three.
First, for many years all the known elements of the kernel of the Magnus
representation arose from considering two separating twists that do
not commute, but whose images under the Magnus representation
do. (We remark that Suzuki has recently constructed other elements of the kernel in \cite{SuzukiNew}.) 
Theorem~\ref{thm:main} says that for separating twists, these are the
only relations between two elements. What happens if we
consider relations between three twists?

\begin{question}
  What words in three positive separating multitwists
  $T_{C},T_{D},T_{E}$ lie in the kernel of the
  Magnus representation $r$?
\end{question}

Although Proposition~\ref{tracecalculation} is still a useful tool\ in
studying this question, it seems that a greater understanding of the
higher intersection form is needed to answer it.

Second, we have focused on the restriction of the Magnus
representation $r$ to the subgroup of the Torelli group generated by
separating twists. But as discussed in the introduction, the Torelli
group has another type of standard generator, the bounding pair
maps.  It would thus be interesting to try to replace one or
both of the separating twists with bounding pair maps:

\begin{question}
  What words in two bounding pair maps lie in the kernel of the
  Magnus representation $r$? In one bounding pair map and one positive separating multitwist?
\end{question}

Finally, one can study other related Magnus representations $r_k$ of
subgroups of the mapping class group, which can be obtained by using
covering spaces other than the universal abelian cover (see Suzuki
\cite{SuzukiGeometric}). Let $\Gamma_k$ be the lower central series of
$\Gamma=\Gamma_1$, defined by $\Gamma_{k+1}=[\Gamma,\Gamma_k]$. The
Johnson filtration consists of the groups $\mathcal{I}_{g,1}(k)$,
where $\mathcal{I}_{g,1}(k)$ is the subgroup of $\MCG$ which acts
trivially on the nilpotent quotient $N_k=\Gamma/\Gamma_k$. Note that
$\mathcal{I}_{g,1}(1)$ is $\MCG$ and $\mathcal{I}_{g,1}(2)$ is
$\mathcal{I}_{g,1}$. Johnson \cite{JohnsonAbel} proved that
$\mathcal{I}_{g,1}(3)$ is actually equal to the Johnson kernel
$\mathcal{K}_{g,1}$, the subgroup of $\MCG$ generated by separating
twists. For each $k$, there is a representation
${r_k\co\mathcal{I}_{g,1}(k)\to \GL_{2g}(\Z[N_k])}$. In particular,
there is a representation of the Johnson kernel
${r_3\co\mathcal{K}_{g,1}\to \GL_{2g}(\Z[\Gamma/\Gamma_3])}$.

\begin{question}
  When is the commutator of two positive separating multitwists
  $[T_{C},T_{D}]$ in $\ker r_3$? In general, what words
  in $T_C$ and $T_D$ lie in
  $\ker r_3$? In $\ker r_k$?
\end{question}\pagebreak

Although we have phrased these questions in the generality of positive
separating multitwists (and one might also consider extensions to
Magnus-multitwists), answers to these questions are not even known for
single separating twists, and it would be logical to begin by
examining this case.

\bibliographystyle{amsplain}
\bibliography{bibliography}

\end{document}